\documentclass{article}
\usepackage{latexsym}
\usepackage{amssymb}
\usepackage{amsfonts}
\usepackage{amsmath}
\usepackage{graphicx}
\usepackage{amsfonts}
\usepackage{epstopdf}
\usepackage{color}
\usepackage{enumerate}
\usepackage[latin1]{inputenc}
\usepackage{mathtools}
\usepackage{amsthm}
\usepackage[latin1]{inputenc}
\usepackage{mathrsfs}
\usepackage{graphicx}
\usepackage[latin1]{inputenc}
\newcommand{\R}{\mathbb{R}}

\newcommand{\GG}{{\cal G}}

\newcommand{\al}{\alpha}

\newtheorem{teo}{Theorem}[section]
\newtheorem{lema}[teo]{Lemma}

\newtheorem{prop}[teo]{Proposition}
\newtheorem{rema}[teo]{Remmark}

\begin{document}
\title{On the strong convergence of multiple ordinary integrals to multiple Stratonovich integrals}
\date{}
\author{Xavier Bardina$^\dagger$ and Carles Rovira\footnote{The authors are supported by the grant PGC2018-097848-B-I00.}}

\maketitle

$^\dagger${\rm Departament de Matem\`atiques, Facultat de Ci\`encies,
	Edifici C, Universitat Aut\`onoma de Barcelona, 08193 Bellaterra}.
{\tt Xavier.Bardina@uab.cat}
\newline

$\mbox{ }$\hspace{0.1cm} $^*${\rm Departament de Matem\`atiques i Inform\`atica,
	Universitat de Barcelona, Gran Via 585, 08007 Barcelona}. {\tt carles.rovira@ub.edu}

\begin{abstract}
Given $\{W^{(m)}(t), t \in [0,T]\}_{m \ge 1}$ a sequence of  approximations to a standard Brownian motion $W$ in $[0,T]$ such that $W^{(m)}(t)$ converges almost surely to $W(t)$ we show that, under regular conditions on the approximations,  the multiple ordinary integrals with respect to $dW^{(m)}$ converge to the multiple Stratonovich integral.
We are integrating functions of the type $$f(x_1,\ldots,x_n)=f_1(x_1)\ldots f_n(x_n) I_{\{x_1\le \ldots \le x_n\}},$$
where for each $i \in \{1,\ldots,n\}$, $f_i$ has continuous derivatives in $[0,T].$ We apply this result to approximations obtained from uniform transport processes.
\end{abstract}

\section{Introduction}

Let $W:=\{W(t), t \ge 0\}$ be a standard Brownian motion and let $\{W^{(m)}(t), t \in [0,T]\}_{m \ge 1}$ a sequence of approximations to $W$ in $[0,T]$ such that $W^{(m)}(t)$ converge almost surely to $W(t)$ as $m$ tends to $\infty$. It is well known that if the approximations are of bounded variation and continuous then, almost surely
\begin{eqnarray*}
& & \lim_{m \to \infty} \int_0^t  W^{(m)}(s) dW^{(m)}(s)= \lim_{m \to \infty} \frac{(W^{(m)}(t))^2}{2} = \frac{(W(t))^2}{2}
\\ & & \qquad \qquad =\int_0^t W(s) dW(s) + \frac{t}{2}= \int_0^t W(s) d^\circ W(s),
\end{eqnarray*}
$t \in [0,T]$, where $dW$ denotes the It\^o integral and $d^\circ W$ is the Stratonovich one.

Wong and Zakai \cite{WZ} extended this result obtaining relations between limits of integrals of the type $\int_0^t \psi(W^{(m)}(s),s) dW^{(m)}(s)$ and $\int_0^t \psi(W(s),s) d^\circ W(s)$. More precisely, assuming that  the approximations are of bounded variation, continuous and for each $m$ can be bounded by a finite random variable and that $\psi(\eta,s)$ has continuous partial derivatives with respect to $\eta$ and $s$, they proved that, for $t \in [0,T]$, almost surely,
\begin{eqnarray*}
	& & \lim_{m \to \infty} \int_0^t  \psi(W^{(m)}(s),s) dW^{(m)}(s)
	\\ & & \qquad \qquad =\int_0^t \psi(W(s),s) dW(s) + \frac12 \int_0^t \frac{\partial \psi}{\partial \eta }(W(s),s) ds 	\\ & & \qquad \qquad =\int_0^t \psi(W(s),s) d^\circ W(s).
\end{eqnarray*}

Gorostiza and Griego \cite{art G-G} applied this result to get strong approximations to stochastic integrals using uniform transport processes as approximations to Brownian motion. They also obtain the corresponding rate of convergence.

Our aim is to extend this result to study the problem of the strong convergence, as $m$ tends to $\infty$, of the multiple integral processes
$$
\Big\{ \int_0^t \int_0^{s_{n}} \ldots \int_0^{s_2} f_n (s_n) \ldots f_1(s_1) d W^{(m)} (s_1)\ldots d W^{(m)} (s_n), t \in [0,T] \Big\},$$
where for each $i \in \{1,\ldots,n\}$, $f_i$ has continuous derivatives in $[0,T].$ We will obtain the almost sure convergence to the multiple Stratonovich integral
$$
\Big\{
\int_0^t \int_0^{s_{n}} \ldots \int_0^{s_2} f_n (s_n) \ldots f_1(s_1) d^\circ W(s_1) \ldots d^\circ W(s_n) , t \in [0,T] \Big\}.$$
We will also apply this convergence to approximations obtained from uniform transport processes. This kind of approximations of multiple integrals has been studied by Bardina and Jolis \cite{BJ} where they obtain a weak convergence result.  Multiple stochastic integrals are useful to represent functionals of Wienner proces and they appear in stochastic Taylor expansions.

The processes  usually called  uniform
transport processes can be presented as
\begin{equation}
W^{(m)}(t)=m^\frac12 (-1)^{A}\int_0^{t}(-1)^{N(um)}du,\label{pois}
\end{equation}
where $\{N(t),\,t\geq0\}$ is a standard Poisson process  and $A\sim
\textrm{Bernoulli}\left(\frac12\right)$  independent of the Poisson
process $N$.

There is some literature  related with strong convergence using approximations based on 
realizations of these processes.
Griego, Heath and Ruiz-Moncayo  \cite{art G-H-RM}  showed that these
processes converge strongly and uniformly on bounded time intervals
to Brownian motion. In  \cite{art G-G2} Gorostiza and Griego
extended the result to diffusions. Again Gorostiza and Griego
\cite{art G-G} and Cs\"{o}rg\H{o} and Horv\'ath \cite{CH} obtained a
rate of convergence. 
Garz\'on, Gorostiza and Le\'on \cite{GGL} defined a sequence of
processes that converges strongly to fractional Brownian motion
uniformly on bounded intervals, for any Hurst parameter $H\in(0,1)$. In \cite{GGL2} and \cite{GGL3}
the same authors deal with subfractional Brownian motion and
fractional stochastic differential equations. Garz\'on, Torres and Tudor \cite{GTT}  also studied the
Rosenblatt process. In \cite{BFR} we get strong approximations to Brownian sheet and in \cite{BBR} for complex Brownian motion.

The proof of our results follows the ideas in \cite{art G-H-RM},   using  It\^o's formula as a main tool. 
It is important to get an expression of the multiple Stratonovich integral and the relation with the multiple It\^o integral well-posed for our problem. We give a detailed expression of this relationship, inspired by the well-known  Hu-Meyer formula  and the work of Kloeden and Platen \cite{KP}.

The paper is organized in the following way. Section 2 is devoted to define the processes  and to give the main results. In Section 3
we study the relationship between It\^o integral and Stratonovich integral. In Section
4  we give the proof of the main theorem.

\section{Notations and main result}

Consider
$$f(x_1,\ldots,x_n)=f_1(x_1)\ldots f_n(x_n) I_{\{x_1\le \ldots \le x_n\}},$$
where $f_i \in L^2[0,T]$ for all $i \in \{1,\ldots,n\}$. From Lemma A.1 and A.2 of \cite{BJ} it is known that $f$ is Stratonovich integrable and the process
\begin{eqnarray*}
I_n^S(t)&:=& I_n^S (f_1,\ldots,f_n)(t) \\  &=& \int_0^t \int_0^{s_{n}} \ldots \int_0^{s_2} f_n (s_n) \ldots f_1(s_1) d^\circ W(s_1) \ldots d^\circ W(s_n)
\end{eqnarray*}
has a version with continuous paths. 
Moreover, there exist the following iterated simple integrals
\begin{equation}\label{relBJ}
Y_k(t) = \int_0^t f_k(s) Y_{k-1}(s) d^\circ W(s),
\end{equation}
for $k \in \{2,\ldots,n\},$ where $Y_1(t)=\int_0^t f_1(s) d^\circ W (s)$ and all these integrals have a continuous version and $Y_n$ coincides with 
$I_n^S$. 

\bigskip

Set, for each $m$, $W^{(m)} (t)$  a sequence of approximations to $W(t)$ such that $W^{(m)} (t)$ is of bounded variation, continuous and converges a.s. to $W(t)$ as $m \to \infty$.
Then, we define the ordinary multiple integral
\begin{eqnarray*}
J_n^{(m)}(t) & := & J_n^{(m)} (f_1,\ldots,f_n)(t) \\
& = & \int_0^t \int_0^{s_{n}} \ldots \int_0^{s_2} f_n (s_n) \ldots f_1(s_1) d W^{(m)} (s_1)\ldots d W^{(m)} (s_n).
\end{eqnarray*}
Obviously, we can consider the iterated integrals
\begin{equation}\label{relm}
J_k^{(m)} (f_1,\ldots,f_k)(t) = \int_0^t f_k(s) J_{k-1}^{(m)} (f_1,\ldots,f_{k-1})(s) d W^{(m)} (s),
\end{equation}
for $k \in \{2,\ldots,n\}$ and where $$J_1^{(m)}(f_1) (t)=\int_0^t f_1(s) d W^{(m)} (s).$$

Following the ideas of Wong and Zakai \cite{WZ}, we  state the set of hypothesis (H) on $W^{(m)}$ as follows:
\begin{enumerate}
	
	\item[(H1)]  almost surely, $W^{(m)}$ are continuous and of bounded variation, 
	
	\item[(H2)]   almost surely, $\lim_{m \to \infty} W^{(m)} (s)=W (s)$ for all $s \in [0,T]$,
	
	\item[(H3)]  for almost all $\omega$, there exists $m_0 (\omega)$ and $K(\omega)$, both finite, such that for all $m > m_0$ and all $s \in [0,T], W^{(m)}(s,\omega) \le K(\omega)$,
	
	\item[(H4)]   $\lim_{m \to \infty} \sup_{s \in [0,T]} \vert W^{(m)} (s)- W (s) \vert = 0$ almost surely,
\end{enumerate}

Then, the main theorem reads as follows:

\begin{teo}\label{mainteo}
	Assume that $f:[0,T]^n \to \R$ is a function
	$$f(x_1,\ldots,x_n)=f_1(x_1)\ldots f_n(x_n) I_{\{x_1\le \ldots \le x_n\}},$$
	such that for each $i \in \{1,\ldots,n\}$, $f_i$ has continuous derivatives in $[0,T].$ 
	Let $(W^{(m)})$ a family of approximations of $W$ that satisfy hypothesis (H1), (H2) and (H3). Then, almost surely,
	$$
	\lim_{m \to \infty} J_n^{(m)}(t) = I^S_n(t)$$
	for all $t \in [0,T]$.
	
	If (H4) is also true, then $$\lim_{m \to \infty} \sup_{s \in [0,T]} \vert J_n^{(m)}(s) -I^S_n(s) \vert = 0$$ almost surely,
\end{teo}

\medskip

We can give now our theorem using approximations based on uniform transport process.

\begin{teo}\label{mainteo2}
	Consider $W^{(m)}:= \{W^{(m)}(t), t \in [0,T]\}$  versions of the uniform transport process (\ref{pois}) on the same probability space as a Brownian motion 
	$\{W(t), t \in [0,T]\}$ such that
	$$\lim_{m \to \infty} \sup_{s \in [0,T]} \vert W^{(m)}(s) -W(s) \vert = 0$$ almost surely.
	Assume that $f:[0,T]^n \to \R$ is a function
	$$f(x_1,\ldots,x_n)=f_1(x_1)\ldots f_n(x_n) I_{\{x_1\le \ldots \le x_n\}},$$
	such that for each $i \in \{1,\ldots,n\}$, $f_i$ has continuous derivatives in $[0,T].$ 
	Then $$\lim_{m \to \infty} \sup_{s \in [0,T]} \vert J_n^{(m)}(s) -I^S_n(s) \vert = 0$$ almost surely,
\end{teo}

\begin{rema}
	Our result can be used to simulate multiple  Stratonovich stochastic integrals $I_n^S$ using as an approximation $ J_n^{(m)}(t)$. In order to compute 
	 $ J_n^{(m)}(t)$ we can use the following decomposition that we proof in Lemma \ref{lem22}:
	 \begin{eqnarray*}
	 	& & J_{n}^{(m)} (t)  =	 \sum_{k=1}^n (-1)^{k+1} J_{n-k}^{(m)} (t) \frac{(W^{(m)}(t))^k}{k!}  (\prod_{l=1}^k f_{n+1-l})(t) \\
	 	& & +  \sum_{k=1}^n (-1)^{k} \int_0^t J_{n-k}^{(m)} (s) \frac{(W^{(m)}(s))^k}{k!}   (\prod_{l=1}^k f_{n+1-l})'(s) ds ,
	 \end{eqnarray*}
	 for any $t \in [0,T]$ and with $J_{0}^{(m)} (t) = 1$.
	
\end{rema}

\begin{rema}
In Bardina and Jolis \cite{BJ} the authors prove the weak convergence, when $m$ tends to infinty, of the multiple integral processes $\{J_n^{(m)}(t), t \in [0,T]\}$ in the space of the real continuous functions when $f_i \in L^2[0,T]$ for all $i \in \{1,\ldots,n\}$.
To obtain the almost sure convergence we need to assume stronger conditions on functions $f_i$ in order to be able to use It\^o formula. The derivability hypothesis that we assume is similar to the one used in the paper of Wong and Zakai \cite{WZ}.

Bardina and Jolis \cite{BJ} also study the weak limit for other classes of integrals. Particularly, they consider the limit of  integrals
$$
\Big\{ \int_0^t \int_0^{s_{n}} \ldots \int_0^{s_2} f(s_1,\ldots,s_n)  d W^{(m)} (s_1)\ldots d W^{(m)} (s_n), t \in [0,T] \Big\},$$
when $f$ is given by a multimeasure. Notice that in this case the almost sure convergence is obtained easily from Theorem 3.1 in \cite{BJ} and the almost sure convergence  in \cite{art G-H-RM}.
\end{rema}

\section{Multiple It\^o and Stratonovich integrals}

The relation between simple It\^o and Stratonovich integrals have been studied deeply. For instance, it is well known that if  $X:=\{X(t), t \ge 0\}$ is an adapted continuous process then
\begin{equation*}
\int_0^t X(s) d^\circ W(s) = \int_0^t X(s) dW(s) + \frac12 [X,W]_t.
\end{equation*}
where $[\,,\,]$ denotes the quadratic covariation.
Particularly, if
$$X(t)=X(0)+\int_0^t u(s) dW(s) + \int_0^t v(s) ds$$
where $\{u(t), t \ge 0\}$ and $\{v(t), t \ge 0\}$ are adapted processes,
and $g$ is a differentiable real function, since we can write that
$$
g(t)X(t)=g(0)X(0)+\int_0^t g(s) dX(s)+\int_0^t X(s) g'(s) ds,$$
we obtain that
\begin{equation}
\int_0^t g (s)X(s) d^\circ W(s) = \int_0^t g(s) X(s) dW(s) + \frac12 \int_0^t g(s) u(s) ds.
\label{itostrato}
\end{equation}

\bigskip
There exists also literature about the relations for  the multiple stochastic integral, begining with the well-known Hu-Meyer formula. Here we obtain some results well adapted to our problem following the notation presented by Kloeden and Platen \cite{KP}.

\smallskip

Let us introduce some notation in order to deal with multiple integrals.
We assume again that $f_i \in L^2[0,T]$ for all $i \in \{1,\ldots,n\}$.
Set
$$
\GG_n:=\big\{ (\al_1,\ldots,\al_m) \, {\rm such \, that} \, \al_i \in \{1,2\}\, \forall i, \, \sum_{i=1}^m \al_i = n \big\}.$$
Notice that if $\alpha=(\alpha_1,\ldots,\alpha_m) \in \GG_n $ then
$\frac{n}{2} \le m \le n$.
Then for any $\alpha=(\alpha_1,\ldots,\alpha_m) \in \GG_n $, let us define
\begin{equation}
I_{(\alpha_1,\ldots,\alpha_m)} (f_1,\ldots,f_n) (t) = \left\{
\begin{array}{lcl}
\int_0^t f_n(s) I_{(\alpha_1,\ldots,\alpha_{m-1})} (f_1,\ldots,f_{n-1})(s) dW(s) \nonumber \\  \qquad \qquad \qquad \qquad \qquad {\rm if} \quad  \alpha_m=1, \nonumber\\ 
\\ \int_0^t f_n(s) f_{n-1}(s) I_{(\alpha_1,\ldots,\alpha_{m-1})} (f_1,\ldots,f_{n-2})(s) ds \nonumber  \\ \qquad \qquad \qquad \qquad \qquad {\rm if} \quad  \alpha_m=2, \label{eqdef}
\end{array}
\right.\end{equation}
with
\begin{eqnarray*}
	I_{(1)}(f_1)(t) &=& \int_0^t f_1(s) dW(s), \\
	I_{(2)}(f_1,f_2)(t) &=& \int_0^t f_2(s) f_1(s) ds.
\end{eqnarray*}	
Clearly, under our assumptions, all these integrals are well-defined. Moreover, for $\alpha=(1,\ldots,1)$ we have the iterated It\^o integral that coincides with the classical multiple It\^o integral:
$$
I_{(1,\ldots,1)} (f_1,\ldots,f_n)(t)=\int_0^t \int_0^{s_{n}} \ldots \int_0^{s_2} f_n (s_n) \ldots f_1(s_1) dW(s_1) \ldots dW(s_n).$$
We will also use the notation 
$$
I_n(t)= I_{n} (f_1,\ldots,f_n)(t) = I_{(1,\ldots,1)} (f_1,\ldots,f_n)(t).$$

We finish this section with three propositions that give us three differents ways to express the Stratonovich integral
$I_n^S(t)= I_{n}^S (f_1,\ldots,f_{n})(t)$. In the first  one we show an expression using the notation we have introduced above. 

\begin{prop}\label{propgn}
Assume that $f_i \in L^2[0,T]$ for all $i \in \{1,\ldots,n\}$. 
Then
$$I_n^S (f_1,\ldots,f_n)(t) =\sum_{(\al_1,\ldots,\al_m) \in \GG_n}
\frac{1}{2^{n-m}}I_{(\alpha_1,\ldots,\alpha_m)} (f_1,\ldots,f_n) (t),$$
for any $t \in [0,T].$
\end{prop}

{\bf Proof:}
We will check the equality by induction on $n$. For $n=1$ it is easy to see that
$$
I_1^S(f_1)(t) = \int_0^t f_1(s) d^\circ W(s) = \int_0^t f_1(s) dW(s) =I_{(1)} (f_1)(t)$$
and that $\{\GG_1\} = \{(1)\}.$

Let us assume that it is true until $n$ and we will check wthat happens for $n+1$. Using (\ref{relBJ}) and the induction hypothesis we have that
\begin{eqnarray}
& &I_{n+1}^S (f_1,\ldots,f_n,f_{n+1})(t) = \int_0^t f_{n+1}(s) I_{n}^S (f_1,\ldots,f_{n})(s) d^\circ W(s) \nonumber \\
& & \quad = \int_0^t f_{n+1}(s) \Big( \sum_{(\al_1,\ldots,\al_m) \in \GG_n}
\frac{1}{2^{n-m}}I_{(\alpha_1,\ldots,\alpha_m)} (f_1,\ldots,f_n) (s)   \Big) d^\circ W(s) \nonumber \\
& & \quad =  \sum_{(\al_1,\ldots,\al_m) \in \GG_n}
\frac{1}{2^{n-m}} \int_0^t f_{n+1}(s) I_{(\alpha_1,\ldots,\alpha_m)} (f_1,\ldots,f_n) (s) d^\circ W(s). \label{eqn+1}
\end{eqnarray}

Notice that if $\al_m = 1$, using (\ref{itostrato}), we can write
 \begin{eqnarray}
 & &\int_0^t f_{n+1}(s) I_{(\alpha_1,\ldots,\alpha_m)} (f_1,\ldots,f_n) (s) d^\circ W(s) \nonumber \\
 & & \quad = \int_0^t f_{n+1}(s) I_{(\alpha_1,\ldots,\alpha_m)} (f_1,\ldots,f_n) (s) d W(s) \nonumber \\ 
 & & \quad \quad + \frac12 \int_0^t f_{n+1}(s) f_{n}(s) I_{(\alpha_1,\ldots,\alpha_{m-1})} (f_1,\ldots,f_{n-1}) (s) ds \nonumber 
 \\
 & & \quad = I_{(\alpha_1,\ldots,\alpha_m,1)} (f_1,\ldots,f_n,f_{n+1}) (t)  
\nonumber \\  & & \quad \quad
  + \frac12  I_{(\alpha_1,\ldots,\alpha_{m-1},2)} (f_1,\ldots,f_{n-1},f_n,f_{n+1}) (t). \label{eqal1}  
 \end{eqnarray}
 On the other hand, using (\ref{itostrato}), when $\al_m=2$ it holds that
  \begin{eqnarray}
 & &\int_0^t f_{n+1}(s) I_{(\alpha_1,\ldots,\alpha_m)} (f_1,\ldots,f_n) (s) d^\circ W(s) \nonumber  \\
 & & \quad = \int_0^t f_{n+1}(s) I_{(\alpha_1,\ldots,\alpha_m)} (f_1,\ldots,f_n) (s) d W(s) 
 \nonumber \\
 & & \quad =  I_{(\alpha_1,\ldots,\alpha_m,1)} (f_1,\ldots,f_n,f_{n+1}) (t).  \label{eqal2}
 \end{eqnarray}
 Putting together (\ref{eqn+1}), (\ref{eqal1}) and (\ref{eqal2}) we finish the proof since
 \begin{eqnarray}
 & &I_{n+1}^S (f_1,\ldots,f_n,f_{n+1})(t) \nonumber \\
 & & \quad =  \sum_{(\al_1,\ldots,\al_m) \in \GG_n}
 \frac{1}{2^{n-m}} I_{(\alpha_1,\ldots,\alpha_m,1)} (f_1,\ldots,f_n,f_{n+1}) (t) \nonumber \\
  & & \quad \quad +  \sum_{(\al_1,\ldots,\al_m) \in \GG_n, \al_m=1}
  \frac{1}{2^{n-m}}  \frac12 I_{(\alpha_1,\ldots,\alpha_{m-1},2)} (f_1,\ldots,f_n,f_{n+1}) (t)\nonumber \\
  & & \quad =  \sum_{(\al_1,\ldots,\al_m) \in \GG_{n+1}}
 \frac{1}{2^{n+1-m}} I_{(\alpha_1,\ldots,\alpha_m)} (f_1,\ldots,f_n,f_{n+1}) (t). \nonumber 
 \end{eqnarray}
 Notice that in the last equality we have used that
 \begin{eqnarray*}
 & &\{ (\al_1,\ldots,\al_m,1); (\al_1,\ldots,\al_m) \in \GG_n\} \\
 && \qquad \cup
 \{ (\al_1,\ldots,\al_{m-1},2); (\al_1,\ldots,\al_m) \in \GG_n \, {\rm with} \, \al_m=1 \}  \\
 & & \quad =
 \{ (\al_1,\ldots,\al_{m+1}) \, {\rm with} \,  \al_i \in \{1,2\} \, \forall i, \al_{m+1}=1, \sum_{i=1}^{m+1} \al_i=n+1 \} \\
&& \qquad \cup
 \{ (\al_1,\ldots,\al_{m}) \, {\rm with} \,  \al_i \in \{1,2\} \, \forall i, \al_{m}=2, \sum_{i=1}^{m} \al_i=n+1 \} 
  \\
 & & \quad = \GG_{n+1}.
\end{eqnarray*}
\hfill$\square$
 
Now, applying Proposition \ref{propgn} we get an iterative definition of $I_n^S$ in terms of $I_{n-1}^S$ and $I_{n-2}^S$.
 
 \begin{prop}\label{propgn2}
 	Assume that $f_i \in L^2[0,T]$ for all $i \in \{1,\ldots,n\}$. Then
 	$$
 I_n^S(t)=\int_0^t f_n(s) I_{n-1}^S (s) dW_s + \frac12 
 \int_0^t f_n(s) f_{n-1}(s) I_{n-2}^S (s) ds,$$
 for any $t \in [0,T].$
 \end{prop}
 
 {\bf Proof:}
It follows easily from the fact that
\begin{eqnarray*}
	&& I_n^S (f_1,\ldots,f_n)(t) =\sum_{(\al_1,\ldots,\al_m) \in \GG_n}
\frac{1}{2^{n-m}}I_{(\alpha_1,\ldots,\alpha_m)} (f_1,\ldots,f_n) (t) \\
&& \, = \sum_{(\al_1,\ldots,\al_m) \in \GG_n, \al_m=1 }
\frac{1}{2^{n-m}} \int_0^t f_n(s) I_{(\alpha_1,\ldots,\alpha_{m-1})} (f_1,\ldots,f_{n-1}) (s) dW(s) \\
&& \,\,  + \sum_{(\al_1,\ldots,\al_m) \in \GG_n, \al_m=2}
\frac{1}{2^{n-m}} \\ && \qquad \times \int_0^t f_n(s) f_{(n-1)} (s) I_{(\alpha_1,\ldots,\alpha_{m-1})} (f_1,\ldots,f_{n-2}) (s) ds
\\
&& \, =  \int_0^t f_n(s) \\
&& \qquad \times\Big(\sum_{(\al_1,\ldots,\al_{m-1}) \in \GG_{n-1}}
\frac{1}{2^{n-1-(m-1)}} I_{(\alpha_1,\ldots,\alpha_{m-1})} (f_1,\ldots,f_{n-1}) (s) \Big) dW(s) \\
&& \,\,  + \frac12  \int_0^t f_n(s) f_{(n-1)} (s) \\ 
\\
&& \qquad \times\Big( \sum_{(\al_1,\ldots,\al_{m-1}) \in \GG_{n-2},}
\frac{1}{2^{n-2-(m-1)}}   I_{(\alpha_1,\ldots,\alpha_{m-1})} (f_1,\ldots,f_{n-2}) (s) \Big) ds.
\end{eqnarray*}
\hfill$\square$

In the next Proposition we express $I_n^S$ using $I_k^S$ for all $k \in \{1,\ldots,n-1\}$. This Proposition is inspired by the work of Wong and Zakai \cite{WZ} and it is the key of the proof of our main theorem.

\begin{prop}\label{lem21}
Assume for each $i \in \{1,\ldots,n\}$, $f_i$ has continuous derivatives in $[0,T].$ Then
\begin{eqnarray*}
	& &	I_n^S(t)= \sum_{k=1}^n (-1)^{k+1} I_{n-k}^{S} (t) \frac{W (t)^k}{k!} (\prod_{l=1}^k f_{n+1-l})(t) \\
	& & \quad\quad-  \sum_{k=1}^n (-1)^k  \int_0^t I_{n-k}^{S} (s) \frac{W (s)^k}{k!} (\prod_{l=1}^k f_{n+1-l})'(s) ds,
\end{eqnarray*}
for any $t \in [0,T]$ and where $I^S_0=1$.
	\end{prop}

{\bf Proof:}
Let us consider the function $F(x,y,u)= xyf_n(u).$ 
By It\^o's formula we can write
\begin{eqnarray}
	& &	F(I_{n-1}^{S} (t), W (t), t) = \int_0^t W (s) f_n(s) d
	I_{n-1}^{S} (s) \nonumber \\
	& & \quad \quad  +  \int_0^t I_{n-1}^{S} (s)  f_n(s) d W (s)
	+  \int_0^t I_{n-1}^{S} (s) W (s) f'_n(s) ds \nonumber\\
	& & \quad \quad  +  \frac12  \int_0^t  f_n(s) d[I_{n-1}^{S} (s),  W (s)].\label{cc}
\end{eqnarray}

Let us check first that 
\begin{equation} \int_0^t I_{n-1}^{S} (s)  f_n(s) d W (s)
	+  \frac12  \int_0^t  f_n(s) d[I_{n-1}^{S} (s),  W (s)] 
	= I_{n}^{S} (t),\label{Itos}
	\end{equation}
Using Proposition \ref{propgn2} it suffices to check that
$$
 \int_0^t  f_n(s) d[I_{n-1}^{S} (s),  W (s)]=
\int_0^t f_n(s) f_{n-1}(s) I_{n-2}^S (s) ds,$$
but it is an obvious consequence again of Proposition \ref{propgn2}. So it is clear that 
(\ref{Itos}) holds.

Then, from (\ref{cc}) and  (\ref{Itos})  and the definition of $F$ it follows that
\begin{eqnarray}
	& &	I_n^S(t)= I_{n-1}^{S} (t) W (t) f_n(t) - \int_0^t W (s) f_n(s) d
	I_{n-1}^{S} (s) \nonumber\\
	& & \quad \quad   - \int_0^t I_{n-1}^{S} (s) W (s) f'_n(s) ds.\label{dd}
\end{eqnarray}

Let us study now the term
$$
\int_0^t W (s)  f_n(s) d
	I_{n-1}^{S} (s). $$
Actually, we will study the more general term
$$
 H^S_k(t):=\int_0^t \frac{(W (s))^k}{k!}  (\prod_{l=1}^{k} f_{n+1-l})(s)
dI_{n-k}^{S} (s) .
$$
From Proposition \ref{propgn2} it follows that
\begin{eqnarray}
H_k^S(t)	&=& \int_0^t \frac{(W (s))^{k}}{k!} (\prod_{l=1}^{k} f_{n+1-l})(s) d
	I_{n-k}^{S} (s)  \nonumber\\
	&=&  \int_0^t \frac{(W (s))^{k}}{k!} (\prod_{l=1}^{k+1} f_{n+1-l})(s)
	I_{n-k-1}^{S} (s)  d W(s) \nonumber\\
	&& \qquad
	+
	\frac12 \int_0^t\frac{(W (s))^{k}}{k!} (\prod_{l=1}^{k+2} f_{n+1-l})(s)
	I_{n-k-2}^{S} (s) ds. \label{ee}
\end{eqnarray}
Consider now the function $F(x,y,u)= x \frac{y^{k+1}}{(k+1)!} (\prod_{l=1}^{k+1} f_{n+1-l})(u).$ 
By It\^o's formula we can write
\begin{eqnarray}
	& &	F(I_{n-k-1}^{S} (t), W (t), t) = \int_0^t \frac{(W (s))^{k+1}}{(k+1)!} (\prod_{l=1}^{k+1} f_{n+1-l})(s) d
	I_{n-k-1}^{S} (s) \nonumber\\
	& & \quad \quad  +  \int_0^t I_{n-k-1}^{S} (s)  \frac{(W (s))^{k}}{k!}  (\prod_{l=1}^{k+1} f_{n+1-l})(s) d W (s)
	\nonumber\\
	& & \quad \quad  +  \int_0^t I_{n-k-1}^{S} (s)  \frac{(W (s))^{k+1}}{(k+1)!} (\prod_{l=1}^{k+1} f_{n+1-l})'(s) ds\nonumber\\
	& & \quad \quad  +  \frac12  \int_0^t  \frac{(W (s))^{k}}{k!}(\prod_{l=1}^{k+1} f_{n+1-l})(s) d[I_{n-k-1}^{S} (s),  W (s)].\label{ff}
\end{eqnarray}
Putting together (\ref{ee}) and (\ref{ff}) and using that from Proposition \ref{propgn2}
\begin{eqnarray*}
&&	\frac12  \int_0^t  \frac{(W (s))^{k}}{k!}(\prod_{l=1}^{k+1} f_{n+1-l})(s)  d[I_{n-k-1}^{S} (s),  W (s)] \\ &&= \frac12 \int_0^t \frac{W (s)^k}{k !} (\prod_{l=1}^{k+1} f_{n+1-l})(s)   f_{n-k-1} (s)
I_{n-k-2}^{S} (s) ds
\end{eqnarray*}
we get that
\begin{eqnarray}
	& &H^S_k(t)= I_{n-k-1}^{S} (t) \frac{(W (t))^{k+1}}{(k+1)!} (\prod_{l=1}^{k+1} f_{n+1-l})(t)
	\nonumber\\
	& & \quad \quad  -  \int_0^t I_{n-k-1}^{S} (s)  \frac{(W (s))^{k+1}}{(k+1)!} (\prod_{l=1}^{k+1} f_{n+1-l})'(s) ds
	\nonumber\\
	& & \quad \quad  - H^S_{k+1} (t).\label{hh}
\end{eqnarray}
So, from (\ref{dd}) we can write
\begin{eqnarray*}
& &	I_n^S(t)= I_{n-1}^{S} (t) W (t) f_n(t) - \int_0^t I_{n-1}^{S} (s) W (s) f'_n(s) ds - H^S_1(t)\\
& & \quad  = I_{n-1}^{S} (t) W (t) f_n(t) - \int_0^t I_{n-1}^{S} (s) W (s) f'_n(s) ds \\
& & \quad \quad  - I_{n-2}^{S} (t) \frac{(W (t))^2}{2!} f_n(t) f_{n-1}(t)+ \int_0^t I_{n-2}^{S} (s) \frac{(W (s))^2}{2!} (f_nf_{n-1})'(s)ds \\
& & \quad \quad  + H^S_2(t),
\end{eqnarray*}
and the proof finishes iterating  (\ref{hh}) $n-1$ times.
\hfill$\square$

\medskip
\section{Proof of the main result}

We begin with a technical lemma where we obtain an expression for 
$J_n^{(m)} =J_n^{(m)} (f_1,\ldots,f_n)(t)$, similar to the expression given in Proposition \ref{lem21} for the multiple Stratonovich integral. This result is also inspired by the paper  \cite{WZ}.

\begin{lema}\label{lem22}
Assume that for each $i \in \{1,\ldots,n\}$, $f_i$ has continuous derivatives in $[0,T].$ Then
\begin{eqnarray*}
	& & J_{n}^{(m)} (t)  =	 \sum_{k=1}^n (-1)^{k+1} J_{n-k}^{(m)} (t) \frac{(W^{(m)}(t))^k}{k!}  (\prod_{l=1}^k f_{n+1-l})(t) \\
	& & +  \sum_{k=1}^n (-1)^{k} \int_0^t J_{n-k}^{(m)} (s) \frac{(W^{(m)}(s))^k}{k!}   (\prod_{l=1}^k f_{n+1-l})'(s) ds ,
\end{eqnarray*}
for any $t \in [0,T]$ and with $J_{0}^{(m)} (t) = 1$.
\end{lema}

{\bf Proof:} Consider the function $F(x,y,u)= xyf_n(u).$ 
Since the functions are of bounded variation
\begin{eqnarray}
& &	F(J_{n-1}^{(m)} (t), W^{(m)} (t), f_n(t)) = \int_0^t W^{(m)} (s) f_n(s) d
J_{n-1}^{(m)} (s)  \nonumber \\
& & \quad \quad  +  \int_0^t J_{n-1}^{(m)} (s)  f_n(s) d W^{(m)} (s)
+  \int_0^t J_{n-1}^{(m)} (s) W^{(m)} (s) f'_n(s) ds.  \label{bound1}
\end{eqnarray}
From (\ref{relm}),  
 (\ref{bound1}) can be written
\begin{eqnarray*}
	& & J_{n}^{(m)} (t)  =	J_{n-1}^{(m)} (t) W^{(m)} (t) f_n(t) \\
	& & - \int_0^t W^{(m)} (s) f_n(s) d
	J_{n-1}^{(m)} (s) -  \int_0^t J_{n-1}^{(m)} (s) W^{(m)} (s) f'_n(s) ds,
\end{eqnarray*}
Thus, now we need to study 
\begin{equation}
\int_0^t W^{(m)} (s) f_n(s) d
	J_{n-1}^{(m)} (s).\label{intleb}
	\end{equation}
Actually, we will study the more general integral
\begin{equation*}
	H_k(t)	:= \int_0^t \frac{(W^{(m)} (s))^k}{k!} (\prod_{l=1}^k f_{n+1-l} )(s) d
	J_{n-k}^{(m)} (s).
\end{equation*}
Consider the function $F(x,y,u)= x \frac{y^{k+1}}{(k+1)!} (\prod_{l=1}^{k+1} f_{n+1-l} )(u).$ 
Since the functions are of bounded variation we can write
\begin{eqnarray}\nonumber
& &	F(J_{n-k-1}^{(m)} (t), W^{(m)} (t), t) = \int_0^t  \frac{(W^{(m)} (s))^{k+1}}{(k+1)!}(\prod_{l=1}^{k+1} f_{n+1-l} )(s) d
J_{n-k-1}^{(m)} (s) \\\nonumber
& & \quad \quad  +  \int_0^t J_{n-(k+1)}^{(m)} (s) \frac{(W^{(m)} (s))^{k}}{k!}   (\prod_{l=1}^{k+1} f_{n+1-l} )(s) d W^{(m)} (s)
\\
& & \quad \quad  + \int_0^t J_{n-(k+1)}^{(m)} (s)
\frac{(W^{(m)} (s))^{k+1}}{(k+1)!} (\prod_{l=1}^{k+1} f_{n+1-l} )'(s)ds.\label{bound3}
\end{eqnarray}
Then, using the definition of $F$ and  (\ref{bound3}) we have that
\begin{eqnarray*}
	H_k(t)	&= & \int_0^t \frac{(W^{(m)} (s))^k}{k!} (\prod_{l=1}^k f_{n+1-l} )(s) d
	J_{n-k}^{(m)} (s) \\
	&= &   \int_0^t J_{n-(k+1)}^{(m)} (s) \frac{(W^{(m)} (s))^k}{k!} (\prod_{l=1}^{k+1} f_{n+1-l} )(s) 
	d W^{(m)} (s)\\
	&= & 	 J_{n-(k+1)}^{(m)} (t) \frac{(W^{(m)} (t))^{k+1}}{(k+1)!} (\prod_{l=1}^{k+1} f_{n+1-l} )(t)  \\
	& & \quad \quad  - \int_0^t  \frac{(W^{(m)} (s))^{k+1}}{(k+1)!} (\prod_{l=1}^{k+1} f_{n+1-l} )(s)  d
	J_{n-(k+1)}^{(m)}(s) \\
	& & \quad \quad  - \int_0^t  J_{n-(k+1)}^{(m)} (s) \frac{(W^{(m)} (s))^{k+1}}{(k+1)!} (\prod_{l=1}^{k+1} f_{n+1-l} )'(s)ds
	\\
	&= & 	 J_{n-(k+1)}^{(m)} (t) \frac{(W^{(m)} (t))^{k+1}}{(k+1)!} (\prod_{l=1}^{k+1} f_{n+1-l} )(t)  \\
	& & \quad \quad  - H_{k+1}(t) \\
	& & \quad \quad  - \int_0^t  J_{n-(k+1)}^{(m)} (s) \frac{(W^{(m)} (s))^{k+1}}{(k+1)!} (\prod_{l=1}^{k+1} f_{n+1-l} )'(s)ds.
\end{eqnarray*}

Iterating the same argument $n-1$ we finish the proof.
\hfill $\square$
\medskip

We can give now the proof of the main theorem.

\medskip
{\bf Proof of Theorem \ref{mainteo}:} We have to check that, almost surely
$$
\lim_{m \to \infty}  J_n^{(m)} (f_1,\ldots,f_n)(t) = I_n^{S} (f_1,\ldots,f_n)(t) ,$$  for any $t \in [0,T].$
We use an induction argument. Let us introduce our hypotesis of induction:
\begin{itemize}
\item[($\bar H_j$)] For any $l \le j$, almost surely  $$ 
\lim_{m \to \infty}  J_j^{(m)} (f_1,\ldots,f_j)(t) = I_j^{S} (f_1,\ldots,f_j)(t) ,$$ for all $t \in [0,T]$.
Moreover, for almost all $\omega$ there exists $m_j (\omega)$ and $K_j(\omega)$, both finite, such that for all $m > m_j$ and all $s \in [0,T]$ and all $l \le j, J^{(m)}_l (s,\omega) \le K_j(\omega)$.
\end{itemize}

Let us study ($\bar H_1$) first. It is a well-known that almost surely
$$ \lim_{m \to \infty} J_1^{(m)}(t) =\lim_{m \to \infty} \int_0^t  W^{(m)}(s) dW^{(m)}(s)=  \int_0^t W(s) d^\circ W(s) = I_1^{(S)}(t) ,
$$
for all $t \in [0,T]$. On the other hand, from Lemma \ref{lem22} we have that 
$$
J_1^{(m)}(t) =
 W^{(m)}(t) f_{1}(t) -
 \int_0^t W^{(m)}(s) f'_{1}(s) ds.$$
On the other hand, the boundedness of $J_1^{(m)}$ is a consequence of the boundedness of $W^{(m)}$ and the continuity of $f_1$ and $f_1'$. So, ($\bar H_1$) is clearly true.

\smallskip
 
 Consider now $j>1$. Let us assume now that ($\bar H_{j-1}$) is true and we will check that ($\bar H_{j}$) holds. Using Proposition \ref{lem21} and Lemma \ref{lem22} to get expressions of $I_{S}^{(m)}$ and $J_{j}^{(m)}$, in order to prove the almost sure convergence it is enough to check that
		for any $k \in \{1,\ldots,j\}$, almost surely
\begin{eqnarray}
&&\lim_{m \to \infty}  J_{j-k}^{(m)} (t) \frac{(W^{(m)})^k}{k!}  (t) (\prod_{l=1}^k f_{j+1-l})(t) \nonumber \\ && \qquad\qquad \qquad= I_{j-k}^{S} (t) \frac{W (t)^k}{k!} (\prod_{l=1}^k f_{j+1-l})(t), \label{lim1}\\
		&& \lim_{m \to \infty} \int_0^t J_{j-k}^{(m)} (s) \frac{(W^{(m)})^k}{k!}  (s) (\prod_{l=1}^k f_{j+1-l})'(s) ds \nonumber \\ && \qquad\qquad \qquad= \int_0^t I_{j-k}^{S} (s) \frac{W (s)^k}{k!} (\prod_{l=1}^k f_{j+1-l})'(s) ds. \label{lim2}
	\end{eqnarray}
for any $t \in [0,T]$. The limit (\ref{lim1}) is an evident consequence of ($\bar H_{j-1}$) and  (\ref{lim2}) follows easily using again ($\bar H_{j-1}$) and dominated convergence. Notice that the boundedness of $J_{j}^{(m)}$ is obtained from ($\bar H_{j-1}$) and Lemma \ref{lem22}.

When (H4) is also true, from Proposition \ref{lem21} and Lemma \ref{lem22} we can write for $j \in \{1,\ldots,n\}$
\begin{eqnarray*}
 &&\sup_{s \in [0,T]} \vert J_j^{(m)}(s) -I^S_j(s) \vert \\
  && \le \sum_{k=1}^n \sup_{s \in [0,T]} \vert (J_{n-k}^{(m)} (s) \frac{W^{(m)}(s)^k }{k!} -  I_{n-k}^{S} (s) \frac{W(s)^k}{k!} )   (\prod_{l=1}^k f_{n+1-l})(s) \vert \\
 \quad& & +  \sum_{k=1}^n \int_0^T \vert (J_{n-k}^{(m)} (s) \frac{W^{(m)}(s)^k }{k!} -  I_{n-k}^{S} (s) \frac{W(s)^k}{k!} )   (\prod_{l=1}^k f_{n+1-l})'(s) \vert ds
 \\
  && \le \sum_{k=1}^n \sup_{s \in [0,T]} \vert J_{n-k}^{(m)} (s) \vert  \sup_{s \in [0,T]}\vert \frac{W^{(m)}(s)^k }{k!} -  \frac{W(s)^k}{k!} \vert  \sup_{s \in [0,T]}  \vert(\prod_{l=1}^k f_{n+1-l})(s) \vert \\
  \quad && + \sum_{k=1}^n \sup_{s \in [0,T]} \vert (J_{n-k}^{(m)} (s)  -  I_{n-k}^{S} (s) \vert  \sup_{s \in [0,T]}\vert \frac{W(s)^k}{k!} \vert  \sup_{s \in [0,T]} \vert(\prod_{l=1}^k f_{n+1-l})(s) \vert \\
 \quad& & +  \sum_{k=1}^n \int_0^T  \vert J_{n-k}^{(m)} (s) \vert \times \vert \frac{W^{(m)}(s)^k }{k!} -  \frac{W(s)^k}{k!} \vert  \times \vert(\prod_{l=1}^k f_{n+1-l})'(s) \vert ds \\
 \quad& & +  \sum_{k=1}^n \int_0^T \vert (J_{n-k}^{(m)} (s)  -  I_{n-k}^{S} (s) \vert \times \vert \frac{W(s)^k}{k!} \vert \times  \vert(\prod_{l=1}^k f_{n+1-l})'(s) \vert ds.
\end{eqnarray*}
Using the equality
$( x^{k}-y^{k})=(x-y) \sum_{k=1}^{n} x^{n-k}y^{k-1}$ we have that
$$\vert W^{(m)}(s)^{k}-W(s)^{k} \vert \le \vert W^{(m)}(s)-W(s) \vert \sum_{k=1}^{n} \vert W^{(m)}(s) \vert^{n-k} \vert W(s) \vert^{k-1}.$$ 
Then using that 
\begin{enumerate}
	\item $\sup_{s \in [0.T]} \vert W(s) \vert < + \infty$ almost surely, 
	\item $\sup_{s \in [0,T]}   \big( \vert(\prod_{l=1}^k f_{n+1-l})(s) \vert +  \vert(\prod_{l=1}^k f_{n+1-l})'(s) \vert  \big) < \infty$
	\item for almost all $\omega$ there exists $m_n (\omega)$ and $K_n(\omega)$, both finite, such that for all $m > m_n$ and all $s \in [0,T]$ and all $l \le n, J^{(m)}_l (s,\omega) \le K_j(\omega)$.
	\item $\lim_{m \to \infty} \sup_{s \in [0,T]} \vert W^{(m)}(s) -W(s) \vert = 0$ almost surely
\end{enumerate}
and doing an interation procedure we can finish easily the proof of the theorem.
\hfill$\square$

\medskip

{\bf Proof of Theorem \ref{mainteo2}:} The proof follows easily from Theorem \ref{mainteo} and the strong convergence result in \cite{art G-H-RM}.

\hfill$\square$

\end{document}